\documentstyle[12pt]{amsart}
\newcommand{\g}{\goth}

\newcommand{\gtsl}{\mbox{\g sl}}

\newcommand{\hgtsl}{\mbox{$\hat{\gtsl}$}}

\newcommand{\nc}{\mbox{${\Bbb C}$}}
\newcommand{\nz}{\mbox{${\Bbb Z}$}}

\newcommand{\cF}{\mbox{${\cal F}$}}

\newcommand{\cP}{\mbox{${\cal P}$}}

\newcommand{\vep}{\varepsilon}
\newcommand{\vphi}{\varphi}
\theoremstyle{plain}
 \newtheorem{thm}{Theorem}
 \newtheorem{prop}[thm]{Proposition}
 \newtheorem{lemma}[thm]{Lemma}
 
\theoremstyle{definition}
 \newtheorem{defn}[thm]{Definition}

\theoremstyle{remark}

\begin{document}
\title{ Zeros and poles of quantum current operators and 
the condition of quantum integrability\\} 
\author{Jintai Ding}
\author{Tetsuji Miwa}
\address{Jintai Ding, RIMS, Kyoto University}
\address{Tetsuji Miwa, RIMS, Kyoto University}
\maketitle
\begin{abstract}
For the current realization of the affine quantum groups, a 
simple comultiplication for the quantum current operators was given 
by Drinfeld. With this 
comultiplication, we study the zeros and poles of 
the quantum current operators and present a condition of integrability
on the quantum current of $U_q\left(\hat{\frak sl}(2)\right)$, which is a 
deformation of the corresponding condition for $\hat{\frak sl}(2)$. We also
present the results about the zeros and poles 
of the quantum current operators of $U_q\left(\hat{\frak sl}(n)\right)$.

\end{abstract}
\pagestyle{plain}
\section{Introduction.} 

For any integrable highest weight module of $\hat{\frak sl}(2)$ of
level $k$, the current operators  $e(z)$ and $f(z)$ satisfy the following 
relations:
$$e(z)^{k+1}=f(z)^{k+1}=0,$$
which we call the condition of integrability \cite{LP}. 

Quantum group was discovered by Drinfeld \cite{Dr1} 
and Jimbo\cite{J1} as a new structure in both mathematics and physics. 
 The definition of a  quantum group is given by  the basic generators and 
the relations based on the data coming from the  
corresponding Cartan matrix.
 However for the case of quantum affine algebras, Drinfeld 
presented a different formulation of affine quantum groups with generators 
in the form of current operators\cite{Dr2}, which, for the case of
  $U_q\left(\hat{\frak sl}(2)\right)$,   give us the 
 quantized current operators corresponding to $e(z)$ and $f(z)$ of 
 $\hat{\frak sl}(2)$. One natural problem is to find out 
if it is possible to find a condition of integrability for the quantum current 
operators, which is a deformation of the classical condition of integrability
above. The non-commutativity of those 
quantum current operators makes the problem much more  difficult than 
the classical case. To solve this problem, we need to use 
a new comultiplication given 
by Drinfeld, which we call Drinfeld comultiplication. 

For the new formulation of affine quantum group, Drinfeld 
proposed another 
comulptiplication formula\cite{DF}\cite{DI} based  on such a formulation.
 The fundamental feature of this comultiplication is its simplicity, 
while the comultiplication formula induced from the conventional
 comultiplication 
 can not be written in a closed form with those current operators.
With this comultiplication, we are able to 
study the zeros and poles of quantum current operators for 
integrable modules.  Our main results are Theorem \ref{thm1}
and \ref{thm2} which state that on any level $k$ integrable module of 
$U_q\left(\hat{\frak sl}(2)\right)$ the matrix coefficients of 
$x^+(z_1)x^+(z_2)....x^+(z_{k+1})$ have
zero at $z_2/z_1=z_3/z_2=\ldots=z_{k+1}/z_k=q^2$,
and those of $x^-(z_1)x^-(z_2)....x^-(z_{k+1})$
have zero at $z_1/z_2=z_2/z_3=\ldots=z_k/z_{k+1}=q^2$,
where $x^+(z)$ and $x^-(z)$ are the quantized 
current operators of $U_q\left(\hat{\frak sl}(2)\right)$
corresponding to $e(z)$ and $f(z)$ of $\hat{\frak sl}(2)$, respectively. 
We first study the case of the fundamental 
representations of level 1 for $U_q\left(\hat{\frak sl}(2)\right)$.
Frenkel and Jing used bosonized vertex operators to construct explicit 
realizations of those representations.
It is clear that all the integrable modules can be 
derived from certain tensors of fundamental representations. 
Therefore, with the Drinfeld comultiplication, we can derive the poles and 
zeros of the quantum currents, which naturally leads to 
a condition integrabililty
for the quantum currents of $U_q\left(\hat{\frak sl}(2)\right)$. 
At the end, we present the corresponding results for 
$U_q\left(\hat{\frak sl}(n)\right)$.

\section{}
For the case of affine quantum  groups, Drinfeld gave a realization of
those algebras in terms of operators in the form of current\cite{Dr2}. 
We will first present such a realization for the case of 
$U_q\left(\hat{\frak sl}(n)\right)$.

Let $A=(a_{ij})$ be the Cartan matrix of type $A_{n-1}$.
\begin{defn}
The algebra $U_q(\hgtsl_n)$ is an associative algebra with unit 
1 and the generators: $\vphi_i(-m)$,$\psi_i(m)$, $x^{\pm}_i(l)$, for 
$i=i,...,n-1$, $l\in \nz $ and $m\in \nz_{\geq 0}$ and a central
 element $c$. Let $z$ be a formal variable and 
 $x_i^{\pm}(z)=\sum_{l\in \nz}x_i^{\pm}(l)z^{-l}$, 
$\vphi_i(z)=\sum_{m\in -\nz_{\geq 0}}\vphi_i(m)z^{-m}$ and 
$\psi_i(z)=\sum_{m\in \nz_{\geq 0}}\psi_i(m)z^{-m}$. In terms of the 
formal variables, 
the defining relations are 
\begin{align*}
& \vphi_i(z)\vphi_j(w)=\vphi_j(w)\vphi_i(z), \\
& \psi_i(z)\psi_j(w)=\psi_j(w)\psi_i(z), \\
& \vphi_i(z)\psi_j(w)\vphi_i(z)^{-1}\psi_j(w)^{-1}=
  \frac{g_{ij}(z/wq^{-c})}{g_{ij}(z/wq^{c})}, \\
& \vphi_i(z)x_j^{\pm}(w)\vphi_i(z)^{-1}=
  g_{ij}(z/wq^{\mp \frac{1}{2}c})^{\pm1}x_j^{\pm}(w), \\
& \psi_i(z)x_j^{\pm}(w)\psi_i(z)^{-1}=
  g_{ij}(w/zq^{\mp \frac{1}{2}c})^{\mp1}x_j^{\pm}(w), \\
& [x_i^+(z),x_j^-(w)]=\frac{\delta_{i,j}}{q-q^{-1}}
  \left\{ \delta(z/wq^{-c})\psi_i(wq^{\frac{1}{2}c})-
          \delta(z/wq^{c})\vphi_i(zq^{\frac{1}{2}c}) \right\}, \\
& (z-q^{\pm a_{ij}}w)x_i^{\pm}(z)x_j^{\pm}(w)=
  (q^{\pm a_{ij}}z-w)x_j^{\pm}(w)x_i^{\pm}(z), \\
& [x_i^{\pm}(z),x_j^{\pm}(w)]=0 \quad \text{ for $a_{ij}=0$}, \\
& x_i^{\pm}(z_1)x_i^{\pm}(z_2)x_j^{\pm}(w)-(q+q^{-1})x_i^{\pm}(z_1)
  x_j^{\pm}(w)x_i^{\pm}(z_2)+x_j^{\pm}(w)x_i^{\pm}(z_1)x_i^{\pm}(z_2) \\
& +\{ z_1\leftrightarrow z_2\}=0, \quad \text{for $a_{ij}=-1$}
\end{align*}
where
\[ \delta(z)=\sum_{k\in \nz}z^k, \quad
   g_{ij}(z)=\frac{q^{a_{ij}}z-1}{z-q^{a_{ij}}}\quad \text{about $z=0$} \]
\end{defn}

In \cite{Dr3}, Drinfeld only gave the formulation of the algebra.
If we  extend the conventional comultiplication to 
those current operators,
the result would be a very  complicated formula which  can not be written 
in a closed form with only those current operators. 
However, Drinfeld also gave the Hopf algebra structure 
for such a formulation in an unpublished note \cite{DF1}.

\begin{thm}
The algebra $U_q(\hgtsl_n)$ has a Hopf algebra structure, which is  given 
by the following formulae. 

\noindent{\bf Coproduct $\Delta$}
\begin{align*}
\text{}& \quad \Delta(q^{\frac c2})=q^{\frac c 2}\otimes q^{\frac c 2}, \\
\text{}& \quad \Delta(x_i^+(z))=x_i^+(z)\otimes 1+
            \vphi_i(zq^{\frac{c_1}{2}})\otimes x_i^+(zq^{c_1}), \\
\text{}& \quad \Delta(x_i^-(z))=1\otimes x_i^-(z)+
            x_i^-(zq^{c_2})\otimes \psi_i(zq^{\frac{c_2}{2}}), \\
\text{}& \quad \Delta(\vphi_i(z))=
            \vphi_i(zq^{-\frac{c_2}{2}})\otimes\vphi_i(zq^{\frac{c_1}{2}}), \\
\text{}& \quad \Delta(\psi_i(z))=
            \psi_i(zq^{\frac{c_2}{2}})\otimes\psi_i(zq^{-\frac{c_1}{2}}),
\end{align*}
where $c_1$ means the action of the center on the 
first component and 
$c_2$ means the action of the center on the second component. 

\noindent{\bf Counit $\vep$}
\begin{align*}
\vep(q^c)=1 & \quad \vep(\vphi_i(z))=\vep(\psi_i(z))=1, \\
            & \quad \vep(x_i^{\pm}(z))=0.
\end{align*}
\noindent{\bf Antipode $\quad a$}
\begin{align*}
\text{(0)}& \quad a(q^c)=q^{-c}, \\
\text{(1)}& \quad a(x_i^+(z))=-\vphi_i(zq^{-\frac{c}{2}})^{-1}
                               x_i^+(zq^{-c}), \\
\text{(2)}& \quad a(x_i^-(z))=-x_i^-(zq^{-c})
                               \psi_i(zq^{-\frac{c}{2}})^{-1}, \\
\text{(3)}& \quad a(\vphi_i(z))=\vphi_i(z)^{-1}, \\
\text{(4)}& \quad a(\psi_i(z))=\psi_i(z)^{-1}.
\end{align*}

\end{thm}

We will give the following  example to explain the comultiplication.

\[
\varphi_i(zq^{c_1\over2})\otimes x_i^+(zq^{c_1})
=\sum_{l_1\in-\nz_{\geq0},l_2\in\nz_{\geq0}}
z^{-(l_1+l_2)}q^{-({l_1\over2}+l_2)c}\varphi_i(l_1)\otimes x_i^+(l_2).
\]

It is clear that the comultiplication structure requires certain completion
on the tensor space. For certain  representations, such 
as the $2$-dimensional representations of
$U_q\left(\hat{\frak sl}(2)\right)$
at a special value, 
this comultiplication may  not be  well-defined. 
Nevertheless, for 
any two highest weight representations, this comultiplication
is well-defined, because the action of the operator as a  coefficient  of 
$z^m$ of the currents operators on any element of such 
a module are zero if $m$ is small enough.  
The explicit proof for the theorem above for the case of 
$U_q\left(\hat{\frak sl}(2)\right)$ is given in \cite{DI}.

We will start with the Frenkel-Jing construction of level $1$ representation 
of $U_q\left(\hat{\frak sl}(2)\right)$ on the Fock space. 

Consider an algebra generated by 
$\{ a_{k}|~k\in \nz \setminus \{ 0\} \}$ 
satisfying:
\[ [a_{k},a_{l}]=\delta_{k+l,0}\frac{[2k][k]}{k}. \]
We call it the Heisenberg algebra. 

Let $ \overline{Q}=\nz \alpha $ be the root lattice of ${\frak sl}(2)$. 
Let us define a group algebra $\nc(q) [\overline{\cP}]$, where
 $\overline{\cP}$ is the weight lattice of ${\frak sl}(2)$. 
Let $\overline{\Lambda}_1$ be the fundamental weight of ${\frak sl}(2)$
 and $2\overline{\Lambda}_1=\alpha$. Let $\overline{\Lambda}_0=0$. 

Set
\[ \cF_{i}:=\nc(q) [a_{-k}(~k\in \nz_{>0})]\otimes
            \nc(q) [\overline{Q}]
              e^{\overline{\Lambda}_i}
  . \]
This gives the Fock space. 

The action of operators 
$a_{k},\partial_{\alpha},e^{\alpha}$~$(1\leq j \leq N)$ is given by 
\begin{align*}
a_{k}\cdot f\otimes e^{\beta}& =\begin{cases}
                      a_{k}f \otimes e^{\beta}          & k< 0;  \\
            \text{$[a_{k},f]$}  \otimes e^{\beta} \quad & k> 0,
                                  \end{cases} \\
\partial_{\alpha}\cdot f\otimes e^{\beta}& 
        =(\alpha,\beta)f\otimes e^{\beta}
         \qquad \text{for}~f\otimes e^{\beta}\in \cF_{i}, \\
e^{\alpha}\cdot f\otimes e^{\beta} & 
=f\otimes e^{\alpha+\beta}.
\end{align*}

\begin{thm} The following 
action  on $\cF_{i}$ of 
$U_q\left(\hat {\frak sl}(2)\right)$
gives a  level 1 highest weight  representation
 with  the $i$-th  fundamental
weight. 
\begin{align*}
& \quad x^{\pm}(z)\mapsto
        \exp[\pm \sum_{k>0}\frac{a_{-k}}{[k]}q^{\mp \frac{1}{2}k}z^k]
        \exp[\mp \sum_{k>0}\frac{a_{k}}{[k]}q^{\mp \frac{1}{2}k}z^{-k}]
        e^{\pm \alpha}z^{\pm \partial_{\alpha}+1}, \\
& \quad \vphi(z)\mapsto
        \exp[-(q-q^{-1})\sum_{k>0}a_{-k}z^k]q^{-\partial_{\alpha}}, \\
& \quad \psi(z)\mapsto
        \exp[(q-q^{-1})\sum_{k>0}a_{k}z^{-k}]q^{\partial_{\alpha}}.
\end{align*}
\end{thm}

This implies that on $\cF_{i}$ for the case of 
$U_q(\hat{\frak sl}(2))$
\begin{align}
&x^+(z)x^+(w)= z^2(1-\frac w z)(1-\frac {w}{zq^2}):x^+(z)x^+(w):,\label{xx+}\\
&x^-(z)x^-(w)= z^2(1-\frac w z)(1-\frac {wq^2}{z}):x^-(z)x^-(w):,\label{xx-}\\
&x^+(z)\vphi(w)= 
q^{-2}\frac {1-\frac w{q^{5/2}z}}{1-\frac {q^{3/2}w}z}:\vphi(w)x^+(z):
=\frac {1-\frac w{q^{5/2}z}}{1-\frac {q^{3/2}w}z}\vphi(w)x^+(z),\label{xp}\\
&\psi(w)x^-(z)= 
\frac{1-\frac {q^{5/2}z}w} {1-\frac z{q^{3/2}w}}\psi(w)x^-(z):
=\frac{1-\frac {q^{5/2}z}w} {1-\frac z{q^{3/2}w}}x^-(z)\psi(w).\label{px}
\end{align}

\begin{lemma} 
Set $\c F=\oplus_{i=0,1}\cF_{i}$. Any level $m$ integrable module 
is a submodule of  $\otimes^m \c F$. 
\end{lemma}

It is clear, for the case of $\hat {\frak sl}(2)$, we have that the 
correlation functions of $e(z)e(w)$ and $f(z)f(w)$ have no poles, 
which are always polynomials of $
z,z^{-1},w,w^{-1}$. By the  correlation functions of an operator, we 
mean all the matrix coefficients of the operator. 
For the quantum case, the correlation functions might have poles.
However, the position of poles are restricted. We have

\begin{prop}\label{prop1}
For any level $k\geq1$ integrable module of
$U_q\left(\hat {\frak sl}(2)\right)$, 
the correlation functions of $x^+(z)x^+(w)$ has at most poles at 
$zq^{-2}=w$. 
\end{prop}  

For the proof, set
\[
\Delta^{(0)}=1,\quad \Delta^{(l)}=(\underbrace{1\otimes\cdots\otimes1}_{l-1}
\otimes\Delta)\Delta^{(l-1)}.
\]
From the comultiplication formula,  we know that on  $\otimes^k \c F$ 
we have 
\begin{align*}
&\Delta^{(k-1)}\left(x^+(z)\right)=\sum_{i=1}^kX^{+(k)}_i(z),\\
&X^{+(k)}_i(z)=\vphi(zq^{1/2})\otimes\vphi(zq^{3/2})\otimes\cdots
\otimes\vphi(zq^{i-3/2})\otimes x^+(zq^{i-1})\otimes1\otimes\cdots\otimes1.
\end{align*}

Let $0<i<j\leq k$. The product $X^{+(k)}_i(z)X^{+(k)}_j(w)$ has
a pole at $z=q^2w$ because it contains $x^+(zq^{i-1})\vphi(wq^{i-1/2})$.
The other terms 
do not have any poles.
Thus, with the lemma above, we complete the proof. 

Similarly, we can show that 

\begin{prop}\label{prop2}
For any level $k\geq1$ integrable module of
$U_q\left(\hat{\frak sl}(2)\right)$, 
the correlation functions of $x^-(z)x^-(w)$ has at most poles  at 
$zq^{2}=w$. 
\end{prop}

Now we state our main result.

\begin{thm}\label{thm1}
For any level $k\geq1$ integrable module of
$U_q\left(\hat{\frak sl}(2)\right)$, 
the correlation functions of $x^+(z_1)x^+(z_2)...x^+(z_k)x^+(z_{k+1})$ 
is zero if $z_2/z_1=z_3/z_2=\ldots=z_{k+1}/z_k=q^2$. 
\end{thm}

We will prove that the correlation functions of
\begin{equation}
X^{+(k)}_{a_1}(z_1)\cdots X^{+(k)}_{a_{k+1}}(z_{k+1})
\label{TERM}
\end{equation}
is zero if $z_2/z_1=z_3/z_2=\ldots=z_{k+1}/z_k=q^2$.

Suppose that $a_m<a_{m+1}$ for some $m$.
Then, the $a_m$-th tensor component of (\ref{TERM})
contains $x^+(z_mq^{a_m-1})\vphi(z_{m+1}q^{a_m-1/2})$.
From (\ref{xp}) we see that this product has a zero at $z_{m+1}=q^2z_m$.
Similarly, if $a_m=a_{m+1}$ the product
$x^+(z_mq^{a_m-1})x^+(z_{m+1}q^{a_m-1})$ has a zero at $z_{m+1}=q^2z_m$.
Proposition \ref{prop1} shows that no poles from other terms
cancel these zeros at $z_{m+1}=q^2z_m$.
Therefore, we conclude that (\ref{TERM}) is zero
if $z_2/z_1=z_3/z_2=\ldots=z_{m+1}/z_m=q^2$, unless $a_m>a_{m+1}$
for all $1\leq m\leq k$.
However, because $1\leq a_m\leq k$, the last case never occurs.
Thus we finish the proof. 

Similarly we have
\begin{thm}\label{thm2}
For any level $k\geq1$ integrable module of
$U_q\left(\hat{\frak sl}(2)\right)$,
the correlation functions of $x^-(z_1)x^-(z_2)...x^-(z_k)x^-(z_{k+1})$ 
is zero
if $z_1/z_2=z_2/z_3=\ldots=z_k/z_{k+1}=q^2$.
\end{thm}

For the case of $U_q(\hat {\frak sl}(n))$, we have the following results,
which can be proved as the case of $U_q(\hat {\frak sl}(2))$. 

\begin{prop}\label{prop3}
For any level $k\geq1$ integrable module of
$U_q\left(\hat{\frak sl}(n)\right)$, the correlation functions of 
$x_i^+(z)x_i^+(w)$ has at most poles  at $zq^{-2}=w$; and  the correlation 
functions of $x_{i\pm 1}^+(z)x_i^+(w)$ 
has at most poles  at 
$zq^{}=w$.
the correlation functions of $x_i^-(z)x_i^-(w)$ has at most poles  at 
$zq^{2}=w$; the correlation functions of $x_{i\pm 1}^-(z)x_i^-(w)$ 
has at most poles  at 
$zq^{-1}=w$.
\end{prop}

The proof for this, we need to use Frenkel-Jing construction to construct 
the level 1 representations. We will omit it here, which can be found in 
\cite{FJ}. However we will list the following formulas, which  are the 
key point of the proof. 

On the level 1 representation in \cite{FJ}, we have:

\begin{align*}
&x_i^+(z)x_i^+(w)= z^2(1-\frac w z)(1-\frac {w}{zq^2}):x_i^+(z)x_i^+(w):,
\\
&x_i^-(z)x_i^-(w)= z^2(1-\frac w z)(1-\frac {wq^2}{z}):x_i^-(z)x_i^-(w):,
\\
&x_i^+(z)\vphi_i(w)= 
q^{-2}\frac {1-\frac w{q^{\frac 5 2}z}}{1-\frac {q^{\frac 3 2}w}z}:\vphi_i(w)x_i^+(z):
=\frac {1-\frac w{q^{\frac 5 2}z}}{1-\frac {q^{\frac 3 2}w}z}\vphi_i(w)x_i^+(z),
\\
&\psi_i(w)x_i^-(z)= 
\frac{1-\frac {q^{\frac 5 2}z}w} {1-\frac z{q^{\frac 3 2}w}}:\psi_i(w)x_i^-(z):
=\frac{1-\frac {q^{\frac 5 2}z}w} {1-\frac z{q^{\frac 3 2}w}}x_i^-(z)\psi_i(w).
\end{align*}

If $a_{ij}=-1$,
\begin{align*}
&x_i^\pm(z)x_j^\pm(w)= \frac
{(z-wq^{\pm 1})}{(z- wq^{\mp 1})}:x_i^+(z)x_i^+(w):,\\
&x_i^+(z)\vphi_j(w)= 
\frac {1-\frac w{q^{-\frac 1 2}z}}{1-\frac {q^{-\frac 3 2}w}z}:\vphi_j(w)x_i^+(z):
=q^{-1}\frac {1-\frac w{q^{-\frac 1 2}z}}{1-\frac {q^{-\frac 3 2}w}z}\vphi_j(w)
x_i^+(z),
\\
&\psi_i(z)x_j^-(w)=
\frac{\frac w{zq^{\frac 12}}-1} {\frac w{zq^{\frac {-3} 2}}-1}:
x_j^-(w)\psi_i(z):
=q\frac{\frac w{zq^{\frac 1 2}}-1} {\frac w{zq^{\frac {-3} 2}}-1}x_j^-(w)
\psi_i(z).
\end{align*}

If $a_{ij}=0$,
\begin{align*}
&x_i^\pm(z)x_j^\pm(w)= x_j^\pm(w)x_i^\pm(z)=:x_i^+(z)x_i^+(w):,\\
&x_i^+(z)\vphi_j(w)= \vphi_j(w)x_i^+(z)=:\vphi_i(w)x_i^+(z):
,\\
&\psi_i(z)x_j^-(w)=x_j^-(w)\psi_i(z)=
:\psi_i(w)x_i^-(z):
.
\end{align*}

\begin{thm}\label{thm3}
For any level $k\geq1$ integrable module of
$U_q\left(\hat{\frak sl}(n)\right)$, 
the correlation functions of $x_{a_1}^+(z_1)x_{a_2}^+(z_2)...x_{a_k}^+(z_k)
x_{a_{k+1}}^+(z_{k+1})$ 
is zero if

 (a)$a_i-a_{i+1}=0$ or $\pm 1 $,

(b)$z_{a_i}/z_{a_{i+1}}=q^{}$ for $a_i-a_{i+1}= \pm 1 $,  

(c)$z_{a_i}/z_{a_{i+1}}=q^{-2}$ for $a_i-a_{i+1}=0 $,

(d)$z_{a_i}/z_{a_{j}}\neq q^{-1}$ for $a_i-a_{j}=\pm 1 $ and $i<j$,

(e)$z_{a_i}/z_{a_{j}}\neq q^{2}$ for $a_i=a_j$ and $i<j$ ;

the correlation functions of $x^-_{a_1}(z_1)x_{a_2}^-(z_2)...
x_{a_k}^-(z_k)x_{a_{k+1}}^-(z_{k+1})$ 
is zero if 

(a)$a_i-a_{i+1}=0$ or $\pm 1 $,

(b)$z_{a_i}/z_{a_{i+1}}=q^{-1}$ for $a_i-a_{i+1}= \pm 1 $,  

(c)$z_{a_i}/z_{a_{i+1}}=q^{2}$ for $a_i-a_{i+1}=0$ and

(d) $z_{a_i}/z_{a_{j}}\neq q^{}$ for $a_i-a_{j}=\pm 1 $ and $i<j$,

(e)$z_{a_i}/z_{a_{j}}\neq q^{-2}$ for $a_i=a_j$ and $i<j$.
\end{thm}

In the classical case, the condition of integrability is used \cite{FS} to 
build semi-infinite construction of the corresponding integrable 
representations, we expect that we can use the quantum condition
of integrability to derive similar constructions, which may even help us to 
resolve certain difficulty in the classical case.

\bigskip
{\it Acknowledgement}

\medskip
\noindent
The authors thank B. Feigin and A. Stoyanovsky for discussions.
J.D. is supported by the grant Reward research (A) 08740020 from the 
Ministry of Education of Japan.

\end{document}